\documentclass[reqno,12pt]{amsart}
\usepackage{amsmath, latexsym, amsfonts, amssymb, amsthm, amscd}
\usepackage{mathrsfs,enumerate}

\numberwithin{equation}{section}

\newtheorem{theorem}{Theorem}[section]
\newtheorem{lemma}[theorem]{Lemma}
\newtheorem{prop}[theorem]{Proposition}

\newtheorem{rem}[theorem]{Remark}
\newtheorem{definition}[theorem]{Definition}
\newtheorem{example}[theorem]{Example}

\newtheorem{remark}[theorem]{Remark}

\newcommand{\cF}{{\ensuremath{\mathcal F}} }

\newcommand{\cM}{{\ensuremath{\mathcal M}} }
\newcommand{\cS}{{\ensuremath{\mathcal S}} }

\newcommand{\cX}{{\ensuremath{\mathcal X}} }

\newcommand{\cZ}{{\ensuremath{\mathcal Z}} }

\newcommand{\E}{{\ensuremath{\mathbb E}} }

\newcommand{\N}{{\ensuremath{\mathbb N}} }

\newcommand{\bbP}{{\ensuremath{\mathbb P}} }

\newcommand{\R}{{\ensuremath{\mathbb R}} }

\newfont{\indic}{bbmss12}

\def\un#1{\hbox{{\indic 1}$_{#1}$}}

\newcommand\MI{\operatorname{MI}}
\newcommand\I{\mathcal I}
\renewcommand\H{\operatorname{H}}

\newcommand\eps{\epsilon}
\newcommand\EST{{}}

\begin{document}

\title[Approximate maximizers of intricacy functionals]{Approximate maximizers of intricacy functionals}

\author{J. Buzzi}

\address{Laboratoire de Math\'ematique d'Orsay - C.N.R.S.  (U.M.R. 8628) \& Universit\'e Paris-Sud
\\
Universit\'e Paris-Sud, F-91405 Orsay Cedex, France}
\email{jerome.buzzi\@@math.u-psud.fr}

\author{L. Zambotti}

\address{Laboratoire de Probabilit{\'e}s et Mod\`eles Al\'eatoires (CNRS U.M.R. 7599) and  Universit{\'e} Paris 6
-- Pierre et Marie Curie, U.F.R. Mathematiques, Case 188, 4 place
Jussieu, 75252 Paris cedex 05, France }
\email{lorenzo.zambotti\@@upmc.fr}

\begin{abstract}
G. Edelman, O. Sporns, and G. Tononi introduced in theoretical
biology the \emph{neural complexity} of a family of random variables.
This functional is a special case of
\emph{intricacy}, i.e., an average of the mutual information of subsystems
whose weights have good mathematical properties. Moreover, its maximum value
grows at a definite speed with the size of the system.

In this work, we compute exactly this \emph{speed of growth}
by building "approximate maximizers" subject to an entropy condition.
These approximate maximizers work \emph{simultaneously} for all intricacies.
We also establish some properties of arbitrary approximate maximizers,
in particular the existence of a \emph{threshold} in the size of subsystems
of approximate maximizers:
most smaller subsystems are almost equidistributed, most larger
subsystems determine the full system.

The main ideas are a random construction of almost maximizers with
a high \emph{statistical symmetry} and the consideration of
\emph{entropy profiles}, i.e., the average entropies of sub-systems
of a given size. The latter gives rise to interesting questions of
probability and information theory.
\end{abstract}



\maketitle


\section{Introduction}

\subsection{Neural Complexity, a measure of complexity from theoretical biology}
In \cite{Tononi94}, G. Edelman, O. Sporns and G. Tononi  introduced
the so called {\it neural complexity} of a family of random variables.
It is defined as an average of mutual
information between any subfamily and its complement, see below.
It has been considered from a theoretical and experimental point
of view by a number of authors, see e.g. \cite{Barnett09,Edelman01,Entropy,
Holthausen99,Krichmar04,Seth06,Seth07,Shanahan08,Sporns00,Sporns02,Sporns07,
Tononi94,Tononi96,Tononi99}.

In order to define the neural complexity, we need to recall
two classical definitions. If $X$ is a random variable
taking values in a finite space $E$, then its
{\it entropy} is defined by
\[
\H(X) := -\sum_{x\in E} P_X(x) \, \log(P_X(x)), \qquad
P_X(x):=\bbP(X=x).
\]
Given two random variables defined over the same probability space,
the {\it mutual information} between $X$ and $Y$ is
\[
\MI(X,Y) := \H(X)+\H(Y)-\H(X,Y).
\]
We refer to the Appendix for a review of the main properties of the
entropy and the mutual information.

Edelman, Sporns and Tononi consider systems formed by a finite family
$X=(X_i)_{i\in I}$ and define the following concept of complexity. For any
$S\subset I$, they divide the system into two subsystems:
\[
X_S:=(X_i, i\in S), \qquad X_{S^c}:=(X_i, i\in S^c),
\]
where $S^c:=I\backslash S$. Then they compute the mutual information
$\MI(X_S,X_{S^c})$ and consider the sum
\begin{equation}\label{esT}
\I(X):=\frac1{|I|+1} \sum_{S\subset I} \frac1{\binom{|I|}{|S|}} \,
\MI(X_S,X_{S^c}),
\end{equation}
where $|I|$ denotes the cardinality of $I$. Note that $\I(X)$ is really a function of the \emph{law} of
$X$.

As shown in \cite{buza1}, one can define more general functionals
\[
\I^c(X):= \sum_{S\subset I} c^I_S \, \MI(X_S,X_{S^c}),
\]
which have similar properties,
provided the properties of "exchangeability" and "weak additivity"
still hold, see Sec. \ref{sec:intricacy}. The resulting functionals have been called
{\bf intricacies} in \cite{buza1}.

Using a super-additivity argument, we showed in \cite{buza1} that the maximum
value of any intricacy over systems with a given size grows linearly with
the size. In this paper, we compute exactly this \emph{speed of growth}
by building "approximate maximizers", i.e., families of an increasing number
of random variables taking value in a fixed set and achieving, in the limit,
the maximum intricacy per variable. Moreover, we shall construct in this paper a sequence
of {\it simultaneous} approximate maximizers for all intricacies.

Our construction is probabilistic in a fundamental way. We shall show that
maximizers should approximately satisfy strong symmetries (see Theorem
\ref{thm:threshold}), that cannot be satisfied exactly
(Lemma \ref{lem:no-ideal-pro}). We shall exhibit a {\it random} sequence
of systems, which satisfy such symmetries {\it in law}, and
approximately satisfy the same symmetries almost surely.

\smallbreak

If the family $(X_i)_{i\in I}$ is completely deterministic or, on the contrary,
independent, then every mutual information vanishes and
therefore $\I^\EST(X)=0$. As these examples suggest,
large values of $\I$ require
compromising between randomness and mutual dependence,
i.e., to have non-trivial
correlation between $X_S$ and $X_{S^c}$ for many subsets $S$.
This explains why maximizing this functional is not
a trivial problem.

\subsection{Main Results}
For the sake of simplicity, we state our
results in this introduction only for the neural complexity \eqref{esT},
deferring the analogous results for arbitrary intricacies to Section 5.

First, we need some notations.
The integers $N\geq1$ and $d\geq 2$
will denote respectively the cardinality of the family $(X_i)_{i\in I}$ and of the range of
each $X_i$. Moreover,
 \begin{itemize}
   \item $\Lambda_{N,d}:=\{0,\ldots,d-1\}^N$ is the set of configurations, i.e., of possible
   values for the random vector $X$;
   \item $\cX(d,N)$ is the set of all $\Lambda_{N,d}$-valued random variables $X$, which we
   shall identify with $\cM(d,N)$, set of all probability measures on $\Lambda_{N,d}$.
 \end{itemize}
In particular, we write indifferently $\H(X)$ and $\H(\mu)$,
as well as $\I(X)$ and $\I(\mu)$.
Of course, entropy and intricacy are in fact functions of the
law $\mu$ of $X$ and not of the (random) values of $X$.

Let us state our main results in the case of the neural complexity:
\begin{theorem}\label{thm:max}
Let $\I^\EST(X)$ be the neural complexity \eqref{esT} of Edelman-Sporns-Tononi.
 \begin{enumerate}
 \item We have for all $\mu\in\cM(d,N)$, setting $x_\mu:=\frac{\H(\mu)}
{N\log d}$,
\begin{equation}\label{Imax}
\frac{\I(\mu)}{N\,\log d}  \leq
x_\mu\left(1-x_\mu\right) \leq \frac 14.
\end{equation}
\item The maximum value of the intricacy at fixed size
\[
\I(d,N):=\max_{X\in\cX(d,N)} \I(X)=\max_{\mu\in\cM(d,N)} \I(\mu)
\]
satisfies:
\[
     \lim_{N\to\infty} \frac{\I(d,N)}N  = \frac{\log d}4.
\]
%
 \item For any $x\in[0,1]$, there exists a sequence $\mu^N\in \cM(d,N)$
 approaching the upper bound of point (1), i.e., satisfying:
 \begin{equation}\label{55}
 \lim_{N\to+\infty} \frac{\H(\mu^N)}{N \log d} = x, \qquad
 \lim_{N\to+\infty} \frac{\I(\mu^N)}{N \log d} = x(1-x).
 \end{equation}
 \end{enumerate}
\end{theorem}

\begin{remark}
{\rm We shall actually prove this theorem for arbitrary intricacies
(see Theorem \ref{thm:main}). More precisely, and perhaps unexpectedly,
we shall build, for each $d\geq2$,
a sequence $\mu^N\in\cM(d,N)$ satisfying, \emph{simultaneously for all intricacies $\I^c$,}
 $$
    \lim_{N\to\infty} \frac{\I^c(\mu^N)}{N} =
    \lim_{N\to\infty} \, \max_{\mu\in\cM(d,N)} \frac{\I^c(\mu)}N,
 $$
 see Remark \ref{cano} below.
}
\end{remark}

\begin{rem}
{\rm
All of the above is new, though numerical experiments by previous
authors \cite[Fig. 1]{Tononi94} had suggested the concavity and the symmetry of the maximal
intricacy given the entropy, but not its quadratic form.
}
\end{rem}

While the upper bound \eqref{Imax} follows from direct computations,
the existence of sequences $(\mu^N)_N$ satistying \eqref{55} is
much less trivial and is the main result of this paper.
As shown in Theorem \ref{thm:threshold} below, such sequences
must exhibit a non-trivial behavior, combining a large amount
of local independence and of non-trivial correlation on a global
level.

The existence of {\bf approximate $x$-maximizers}, i.e., sequences
$\mu^N\in\cM(d,N)$ satisfying \eqref{55}, follows in our approach
from a probabilistic construction: we shall prove that uniform
distributions on appropriately chosen {\it random sparse supports}
will have almost surely the desired properties: see Proposition
\ref{duque} below.

\medbreak
In the course of the proof, we also obtain rather detailed information on
the structure of approximate $x$-maximizers. A key notion is the following one.

\begin{definition}\label{def:profile}
Given $X\in\cX(d,N)$, its
{\bf entropy profile} is the function $h_X:[0,1]\to[0,1]$ such that
$h_X(0)=0$,
 $$
    h_X\left(\frac kN\right) = \frac{1}{\binom{N}{k}}
     \sum_{S\subset I, \, |S|=k}
      \frac{\H(X_S)}{N\log d}, \qquad k\in I:=\{1,\dots,N\}
 $$
and $h_X$ is affine on each interval $\left[\frac{k-1}N,\frac kN\right]$,
$k\in I$.
\end{definition}

\begin{theorem}\label{thm:profile}
For $x\in[0,1]$, let the \emph{ideal profile} be
 $$
   h_x^*(t)=x \wedge t=\min\{x,t\}.
 $$
Then for any sequence $\mu_N\in\cM(d,N)$ of approximate $x$-maximizers,
we have
\[
   \|h_{\mu^N}-h_x^*\|_{\sup}:=\sup_{t\in[0,1]} |h_{\mu^N}(t)-h_x^*(t)| \to 0 \text{ as }N\to\infty.
\]
In particular, for any sequence $\mu^N\in \cM(d,N)$
of approximate maximizers, i.e., such that $\lim_{N\to\infty} \I(\mu^N)/N=\log d/4$,
we have:
 $$
  \lim_{N\to+\infty} \frac{\H(\mu^N)}{N \log d} = 1/2,
  \qquad
  \lim_{N\to+\infty} \|h_{\mu^N}-h_{1/2}^*\|_{\sup}=0.
 $$
\end{theorem}
\noindent
Again, we prove in fact a version of this result for all intricacies,
see Theorem \ref{barry} below.

\medbreak
If $(\mu^N)_N$ is a sequence of approximate $x$-maximizers and
$X^N\in\cX(d,N)$ has law $\mu^N$, we say that
$(X^N)_N$ is also a sequence of approximate $x$-maximizers.

\medbreak

A corollary of the convergence of entropy profiles is the existence of
a threshold in the behavior of typical subsystems  of
approximate $x$-maximizers: if $|S|\leq xN$, then $X^N_S$ is
almost uniform, which corresponds to local independence; if
 $|S|\geq xN$, then the whole family $X^N$ is almost a function
of $X^N_S$, which corresponds to strong global correlation.
Recall that $H(Y \, | \, Z)$ is the conditional
entropy of $Y$ given $Z$, see the Appendix below.

\begin{theorem}\label{thm:threshold}
Let $(X^N)_N$ be an approximate $x$-maximizer. Let $y\in]0,1[$
and set $k_N:=\lfloor yN \rfloor$.
Consider
the $\binom{N}{k_N}$ sub-systems $X_S$ of $X^N$ of size $k_N=|S|$.
For all $\eps>0$, if $N$ is large enough, then
except for at most $\eps\binom{N}{k_N}$ of such subsets $S$, the following
holds:
 \begin{itemize}
  \item if $y\leq x$: $X_S$ is
\emph{almost uniform}: ${\displaystyle 1-\eps\leq\frac{\H(X_S)}{|S|\log d}
\leq 1}$;
 \item if $y\geq x$: $X_S$ \emph{almost determines} the whole system $X^N$:
 ${\displaystyle 0\leq \frac{\H(X^N|X_S)}{N\log d}\leq \eps}$.
 \end{itemize}
\end{theorem}
\noindent
Again, we prove a more general version of this result in
Theorem \ref{ott} below.

\subsection{Strategy of Proof and Organization of the paper}

The main ideas of the proofs of the two theorems are a probabilistic construction of the sequence
maximizers and the consideration of the entropy profiles $h_X$ defined above.
As we indicated, in fact we analyze arbitrary
intricacies generalizing neural complexity.

\medbreak

In section 2 we recall the notion of intricacy as
a family of functionals over finite sets of
discrete random variables satisfying exchangeability and
weak additivity and give simple examples. In section 3
we give upper bounds on the intricacies of arbitrary systems
of given size and entropy. In section 4 we
prove the main results by means of a probabilistic construction
of random approximate maximizers. In section 5 we collect
our results for arbitrary intricacies.
An Appendix contains basic facts from entropy theory for the convenience
of the reader.

\subsection{Further questions}

The bound $x(1-x)$ in \eqref{55} is symmetric with respect to $x=1/2$ and
independent of $d\geq2$. We do not know whether these simple properties,
which extend to arbitrary intricacies (see Theorem \ref{thm:main}),
can be proved directly, e.g.: does there exist a duality operation
in $\cX(N,d)$ exchanging systems with entropy $xN\log d$ and $(1-x)N\log d$
while preserving their intricacy? Can one deduce from a system in $\cX(d,N)$
with entropy $H$ and intricacy $I$ a system in $\cX(d',N)$ with
entropy $(\log d'/\log d)H$ and intricacy $(\log d'/\log d)I$?

\medbreak

This work has focused on properties of systems with size tending
to infinity. Notice that we know very little on the exact maximizers
for fixed size beyond the constraints on their entropy contained
in our main results. Because of the invariance properties of intricacy
(see Lemma \ref{lem:invariance} and the following comment), exact
maximizers are non-unique  but we do not even know if there are only finitely
many of them.

\medbreak

Our construction of approximate maximizers is probabilistic.
Could it be done deterministically? Would the corresponding
algorithms possess a computational complexity related to the
complexity that intricacies are supposed to describe?

\medbreak

Our construction is global but could systems with maximum intricacy
be built by a local approach, i.e., a  "biologically reasonable"
building process, using some type of local rules and/or evolution?
That is, does there exist a "reasonable" self-map $T:\cM(d,N)\to\cM(d,N)$
such that the neural complexity of $T^n(\mu)$ converges to the maximum
as $n\to\infty$ for "many" $\mu\in\cM(d,M)$.

\medbreak

Our work also leads to
interesting probabilistic constructions and questions in
the theory of entropy and information. For instance:

\medbreak
\noindent {\bf Problem.} {\it
Describe the set of functions $h:\{0,\dots,N\}\to\mathbb R$
obtained from picking $X\in\cX(d,N)$ and setting $h(k)$ to
be the average entropy of $X_S$ where $S$ ranges
over the subsets of $\{1,\dots,N\}$ with cardinality $k$.
}
\medbreak

Basic properties of the entropy (recalled in the Appendix)
imply that $h(0)=0$,
$0\leq h(k+1)-h(k)\leq \log d$ for $0\leq k<N$ and
 $$
   h(k+\ell)-h(k)\leq h(j+\ell)-h(j) \text{ for }0\leq j\leq k\leq k+\ell\leq N.
 $$
However we shall show that not all such functions $h$ arise from some
$X\in\cX(d,N)$, see Lemma \ref{lem:no-ideal-pro}. See \cite{maka} for
a closely related question.

\section{Intricacy}\label{sec:intricacy}

\subsection{Definition}
In this paper, a {\bf system} is a finite collection $(X_i)_{i\in I}$ of random variables, each $X_i$
taking value in the same finite set $V$. Without loss of generality, we assume that $V=\{0,\dots,d-1\}$
for all $i\in I$ and some $d\geq2$ ($d$ should be thought of as a convenient normalization) and $I$ is a set of positive integers.
We let $\cM=\bigcup_{d\geq2} \cM(d)=\bigcup_{d\geq 2,N\geq1} \cM(d,N)$
be the set of the corresponding laws, that is, the probability measures on
$\{0,\dots,d-1\}^I$ for each finite subset $I\subset \N^*:=\{1,2,3,\dots\}$.

For $S\subset I$, we denote
\[
X_S:=(X_i, i\in S).
\]
In \cite{buza1}, we defined the following family of functionals over such systems (more precisely: over their laws)
formalizing (and slightly generalizing) the neural complexity
of Edelman-Sporns-Tononi \cite{Tononi94}:

\begin{definition}
A {\bf system of coefficients} is
a collection of numbers
\[
c:=(c_S^I: \, I\subset\subset\N^*, \, S\subset I),
\]
i.e., $I$ ranges over the finite subsets of $\N^*$,
for all $I$ and all $S\subset I$:
 \begin{equation}\label{eq:cond-CIS}
     c_S^I\geq 0, \quad \sum_{S\subset I} c_S^I = 1, \quad \text{ and } c^I_{S^c}=c^I_S
 \end{equation}
where $S^c:=I\setminus S$.
The corresponding {\bf mutual information functional} is $\I^c:\cM\to\mathbb R$ defined by:
 $$
    \I^c(X):=\sum_{S\subset I} c^I_{S} \MI\left(X_S,X_{S^c}\right).
 $$
By convention,
$\MI\left(X_\emptyset,X_I\right)=\MI\left(X_I,X_\emptyset\right)=0$.
An {\bf intricacy}, is a mutual information function
satisfying:
 \begin{enumerate}
  \item {\bf exchangeability} (invariance by permutations): if $I,J\subset\subset\N^*$
  and $\phi:I\to J$ is a bijection, then $\I^c(X)=\I^c(Y)$ for any $X:=(X_i)_{i\in I}$, $Y:=(X_{\phi^{-1}(j)})_{j\in J}$;
  \item {\bf weak additivity}: for any two \emph{independent} sub-systems $(X_i)_{i\in I},(Y_j)_{j\in J}$
  (defined on the same probability space): $\I^c(X,Y)=\I^c(X)+\I^c(Y)$.
 \end{enumerate}
$\I^c$ is {\rm \bf non-null} if  some coefficient $c^I_S$ with $S\notin\{\emptyset, I\}$ is not zero.
\end{definition}

\subsection{Classification of intricacies}

In section 3 of \cite{buza1} the following has been proved

\begin{prop}\label{pro:special}
A mutual information functional $\I^c$ determines its coefficients uniquely and the following equivalences hold:
 \begin{itemize}
  \item $\I^c$ is exchangeable if and only if $c^I_S$ depends only on
  $|I|$ and $|S|$;
  \item an exchangeable $\I^c$ is weakly additive if and only if there exists a random variable
  $W_c\mapsto[0,1]$ such that $W_c$ and $1-W_c$ have the same law and
  \begin{equation}\label{mainc}
   c^N_k=\E\left((1-W_c)^{N-k} \, W_c^{k}\right)=\int_{[0,1]}x^k(1-x)^{N-k}\, \lambda_c(dx),
   \end{equation}
   where $\lambda_c$ is the law of $W_c$.
\item an exchangeable weakly additive $\I^c$ is non-null iff $\lambda_c(]0,1[)>0$, in which
case all coefficients $c^I_S$ are non-zero.
  \end{itemize}
\end{prop}
\noindent
In this paper we consider only non-null intricacies $\I^c$.
\begin{example}\label{ccca}{\rm
The intricacy $\I^{\EST}$ of Edelman-Sporns-Tononi  is defined by the
coefficients:
\begin{equation}\label{ET}
    c^I_S = \frac{1}{|I|+1}\frac{1}{\binom{|I|}{|S|}}
\end{equation}
and it is easy to see that in this case \eqref{mainc} holds with $W_c$ a uniform
variable over $[0,1]$, see Lemma 3.8 in \cite{buza1}.
For $0<p<1$, the {\bf symmetric $p$-intricacy} $\I^p$ is defined by
 $$
     c^I_S = \frac12\left(p^{|S|}(1-p)^{|I\backslash S|}+(1-p)^{|S|}\, p^{|I\backslash S|}\right)
 $$
 and in this case $W_c$ is uniform on $\{p,1-p\}$.
For $p=1/2$, this yields the {\bf uniform intricacy} $\I^U(X)$ with:
 $$
     c^I_{S}=2^{-|I|},
 $$
 and $W_c=1/2$ almost surely. All these functionals are clearly non-null
and exchangeable.
 }
\end{example}

\begin{rem}
{\rm
The global $1/(|I|+1)$ factor in \eqref{ET}
is not present in \cite{Tononi94}, which did not compare systems of different sizes.
However it is necessary in order to have weak additivity.
}
\end{rem}

\subsection{Simple examples}\label{tse}

Let $X_i$ take values in $\{0,\ldots,d-1\}$ for all $i\in I$,
a finite subset of $\N^*$.

\begin{example}\label{gb}
{\rm If the variables $X_i$ are independent then each mutual information
is zero and therefore:
 $
    \I^c(X)=0. \qquad \qquad \square
 $
}
\end{example}
\begin{example}\label{gb2}
{\rm If each $X_i$ is a.s. equal to a constant $c_i$ in $\{0,\ldots,d-1\}$,
then, for any $S\ne\emptyset$, $\H(X_S)=0$. Hence,
 $
   \I^c(X) = 0. \qquad \qquad \square
 $
}
\end{example}
\begin{example}\label{gb3}
{\rm If $X_1$ is uniform on $\{0,\ldots,d-1\}$ and $X_i=X_1$ for
all $i\in I$, then,
for any $S\ne\emptyset$, $\H(X_S)=\log d$ and, if additionally $S^c\ne\emptyset$,
$\H(X_S\bigm|X_{S^c})=0$ so that each mutual information $\MI(X_S;X_{S^c})$ is $\log d$.
Hence,
 $$
   \I^c(X) = \sum_{S\subset I\backslash\{\emptyset, I\}} c_S^I \cdot \log d
=\left(1-c_\emptyset^I-c_I^I\right) \log d\leq \log d. \qquad \qquad \square
 $$
}
\end{example}

\medbreak

Examples \ref{gb} and \ref{gb2} correspond to, respectively,
maximal and minimal total entropy. In these extreme cases $\I^c=0$.
Example \ref{gb3} has positive total entropy and intricacy (if
all $c^I_S$ are non zero). However,
the values of the intricacy grow very slowly with $|I|$ in these
examples: they stay bounded.
We shall see however
how to build systems $(X_i)_{i\in I}\in\cX(d,I)$ which realize much
larger values of $\I^c$, namely of the order of $|I|$.

\medbreak

\subsection{Invariance properties of intricacies}
We have the following obvious invariances of intricacies.

\begin{lemma}\label{lem:invariance}
The intricacies are invariant under the following group actions
on $\cX(d,N)$ for some $N,d\geq1$:
 \begin{enumerate}
  \item the group $\cS_N$ of permutations on $\{1,\dots,N\}$
 acting on $\cX(d,N)$ by: \
 $(\sigma X)_i = X_{\sigma^{-1}(i)}$, $\forall \, i=1,\dots,N$.
  \item the $N$th power $(\cS_d)^N$ of the permutation group on $\{0,\dots,d-1\}$
 acting on $\cX(d,N)$ by: \
 $(\sigma X)_i = \sigma_i\circ X_i$, $\forall \,i=1,\dots,N$.
 \end{enumerate}
\end{lemma}

In particular, for $N,d\geq2$, the maximum of $\I^c$ over
$X\in\cX(d,N)$ cannot be achieved at a single probability
measure on $\Lambda_{d,N}=\{0,\dots,d-1\}^N$. Indeed, if it were the case,
then this measure would be invariant under the group
action (2) above. However, this action is transitive on $\Lambda_{d,N}$.
Therefore the measure would be equidistributed on this set.
Hence the maximizer would be a family of independent variables,
for which the intricacies are all zero. This is
a contradiction whenever $N,d\geq2$.

\section{Upper bounds on intricacies}

In \cite{buza1}, it was proved that $\I^c(X)<N\log d/2$ if
$X\in\cX(d,N)$. By comparison with "ideal entropy profiles" defined
below, we prove sharper upper bounds for systems with given size
and entropy.

\subsection{Definitions}
We define the ideal entropy profile and the corresponding
intricacy values both for finite size and in the limit $N\to\infty$.
We also introduce an adapted norm to measure the distance
between profiles.

Let $\I^c$ be some intricacy. It is convenient to use the
following probabilistic representation of the coefficients $c$
based on the the random variable $W_c$ with law $\lambda_c$
defined by \eqref{mainc}. Let $(Y_i)_{i\geq 1}$ be
a sequence of i.i.d. uniform random variables on $[0,1]$ and let
\begin{equation}\label{lava0}
D_{N}:=\sum_{k=1}^N \un{(Y_k\leq W_c)}, \qquad \beta_N:=\frac{D_N}N, \qquad N\geq 1.
\end{equation}
Conditionally on $W_c$, $D_N$ is a binomial variable with parameters $(N,W_c)$. In particular, for all $g:\N\mapsto\R$, by \eqref{mainc}
\begin{equation}\label{lava}
\E\left(g\left(D_N\right)\right) = \int_{[0,1]}\sum_{k=0}^N \binom Nk
\, x^k(1-x)^{N-k} \, g(k) \, \lambda_c(dx)
= \sum_{k=0}^N c^N_k \, \binom{N}{k} \, g(k),
\end{equation}
and therefore, for all bounded Borel $f:[0,1]\mapsto\R$
\begin{equation}\label{lava2}
\E\left(f\left(\beta_N\right)\right) =
\sum_{k=0}^N c^N_k \, \binom{N}{k} \, f\left(\frac{k}N\right).
\end{equation}

We recall the Definition \ref{def:profile} of the entropy profile of
$X\in\cX(d,N)$: $h_X(0)=0$,
 $$
    h_X\left(\frac kN\right) = \frac{1}{\binom{N}{k}}
     \sum_{S\subset I,\, |S|=k}
      \frac{\H(X_S)}{N\log d}, \qquad k\in I:=\{1,\dots,N\}
 $$
and $h_X$ is affine on each interval $\left[\frac{k-1}N,\frac kN\right]$,
$k\in I$. We can now define the ideal profiles and their intricacies.

\begin{definition}
For $x\in[0,1]$ and $N\geq1$,
the {\bf ideal entropy profile} is
 \begin{equation}\label{eq:ideal}
   h^*_x(t):=t\wedge x=\min\{t,x\}
 \end{equation}
and the corresponding (normalized) intricacies are, for finite $N$:
 \begin{equation}\label{eq:max-intric2}
    i^c_N(x):=2\sum_{k=0}^N c^N_k \, \binom{N}{k} \,
  h^*_x(k/N) -x=2\, \E\left(x\wedge \beta_N\right)-x
  \end{equation}
and, for $N\to\infty$:
 \begin{equation}\label{eq:max-intric}
    i^c(x):=2\int_0^1 (t\wedge x) \, \lambda_c(dt)-x=2\, \E(x\wedge W_c)-x.
  \end{equation}
\end{definition}
We remark that the ideal profile $h^*_x$ does not depend on the intricacy $\I^c$.
Finally, we define a family of norms. For all bounded Borel $f:[0,1]\mapsto\R$,
let
 \begin{equation}\label{norm}
   \|f\|_{c,N} := \sum_{k=0}^N c^N_k \, \binom{N}{k}
    \,  |f(k/N)| = \E\left(|f(\beta_N)|\right).
  \end{equation}
\begin{rem}\label{nonni}{\rm
For the particular cases of Example \ref{ccca} we have more explicit expressions.
For the Edelman-Sporns-Tononi neural complexity, the above reduces to
\[
    i(x)=x(1-x), \quad x\in[0,1],
\]
  for the uniform intricacy
\[
    i^U(x)=\min\{x,1-x\}, \quad x\in[0,1],
\]
  and for the symmetric $p$-intricacy
\[
 i^p(x)=\min\{x,1-x,p,1-p\}, \quad x\in[0,1].
\]

}
\end{rem}

\subsection{Upper bounds and distance from the ideal profile}
In this section we prove the following upper bounds

\begin{prop}\label{violi}
Let $\I^c$ be an intricacy.
\begin{enumerate}
  \item $i^c:[0,1]\to[0,1]$ is a concave function admitting the Lipschitz constant $1$
  and symmetric about $1/2$: $i^c(1-x)=i^c(x)$. Moreover, $i^c(1/2)=\max_{x\in[0,1]} i^c(x)$.
  \item $|i^c(x)-i^c_N(x)|\leq 1/\sqrt N$.
  \item All systems $X\in\cX(d,N)$ with $\frac{\H(X)}{N\log d}=x$
   satisfy:
\begin{equation}\label{mar}
\frac{\I^c(X)}{N\log d} = i^c_N(x)-\|h_X-h_x^*\|_{c,N}\leq i^c_N(x).
\end{equation}
 \item If $X^N\in\cX(d,N)$ and
$\lim_{N\to\infty}\frac{\H(X^N)}{N\log d}=x$, then
\begin{equation}\label{mari}
     \limsup_{N\to\infty} \frac{\I^c(X^N)}{N\log d}
        \leq i^c(x).
\end{equation}
\end{enumerate}
\end{prop}

Observe that to show that
$i^c(x)$ is indeed the value of the limit \eqref{mari}, rather than a mere
upper bound, requires to prove the existence of sequences saturating the
inequality.  This is deferred to the
next section. Before proving Proposition \ref{violi}, we need some
preliminary material which will also be useful later.

\subsection{The functions $i^c_N(x)$ and $i^c(x)$}
We consider the first two points of the proposition,
beginning with the convergence of $i^c_N(x)\to i^c(x)$.
\begin{lemma}\label{opti-I}
For all $x\in[0,1]$
\[
|i^c_N(x)-i^c(x)|\leq \frac1{2\sqrt{N}}, \qquad N\geq 1.
\]
\end{lemma}
\begin{proof}
We use the probabilistic representations \eqref{eq:max-intric} and \eqref{eq:max-intric2} and
we obtain
\[
|i^c_N(x)-i^c(x)|\leq \E\left( \left|h^*_x(\beta_N)-h^*_x(W_c)\right|\right)\leq
\E\left( \left|\beta_N-W_c\right|\right)\leq \sqrt{\E\left( \left|\beta_N-W_c\right|^2\right)}.
\]
Since $D_N=N\beta_N$ is, conditionally on $W_c$, a binomial variable with parameters $(N,W_c)$,
we have that
\begin{equation}\label{oo}
\E\left( \left|\beta_N-W_c\right|^2\right)=\E\left( {\rm Var}\left(\beta_N \, |\, W_c\right)\right)=
\E\left( \frac{W_c(1-W_c)}N\right) \leq \frac1{4N}
\end{equation}
and the result is proven.
\end{proof}

\medbreak

Now, we analyze the limit function $i^c(x)$.

\begin{lemma}\label{pro}$ $
\begin{enumerate}
\item $i^c(x)=\E(\min\{x,1-x,W_c,1-W_c\})$ for all $x\in[0,1]$.
 \item The function $i^c:[0,1]\mapsto[0,1]$  is
$1$-Lipschitz and concave. The distributional second derivative of $i^c$
is $-2\lambda_c$.
\item $i^c(x)=i^c(1-x)$ for all $x\in[0,1]$.
\item $i^c$ achieves its maximum at $x=1/2$ and $i^c(1/2)=\E(W_c\wedge (1-W_c))$.
\item $i^c$ is maximum only at $x=1/2$ if and only if $1/2$ belongs to the support of $\lambda_c$.
\end{enumerate}
\end{lemma}

\begin{proof}
First, for all $x,a\in[0,1]$
\begin{equation}\label{ny}
x\wedge a+x\wedge(1-a)-x = \min\{x,1-x,a,1-a\}.
\end{equation}
Indeed, one can assume $a\leq1-a$ and then check the above in
the three cases: $x\leq a$, $a\leq x\leq 1-a$ and $x\geq1-a$.
Since $W_c$ has same law as $1-W_c$, point (1) and (3) follow.

Concavity, $1$-Lipschitz continuity and symmetry w.r.t. $1/2$ follow easily.
Moreover, an integration
by parts shows that for all $\varphi\in C^\infty(\mathbb R)$
with compact support contained in $(0,1)$:
\[
\begin{split}
 & \int_{[0,1]}\varphi''(x)\, i^c(x)\, dx = 2\int_{[0,1]}
\left[\int_{[0,1]}\varphi''(x)\, x\wedge t\ dx \right] \lambda_c(dt)
-\int_{[0,1]}\varphi''(x)\, x\, dx\\
& = 2\int_{[0,1]}
\left[\int_0^t\varphi''(x)\, x\, dx + \int_t^1\varphi''(x)\, t dx \right] \lambda_c(dt)
- 0
\\ & = 2\int_{[0,1]}
\left[\varphi'(t)t-\varphi(t) - t\varphi'(t) \right] \lambda_c(dt)
  =  -2\int_{[0,1]}\varphi(t)\, \lambda_c(dt),
\end{split}
\]
proving that $(d/dx)^2i^c=-2\lambda_c$ as distributions. Point (2) is proved.

Since $i^c$ is concave and $i^c(x)=i^c(1-x)$ then
$i^c(x)\leq i^c(1/2)=\E(W_c\wedge(1-W_c))$ for all $x\in[0,1]$.
Point (4) is proved.

Let us now assume $x<1/2$ (the other case being similar) so that $x=x\wedge(1-x)$. Set $w=W_c\wedge(1-W_c)$. Then, by \eqref{ny}
\[
i^c(1/2)-i^c(x)=\E\left(w-x\wedge w\right)
 =\E\left((w-x)\,
\un{x<w)}\right).
\]
Hence,
$i^c(x)<i^c(1/2)$ if and only if $\bbP(x<w)>0$, i.e.
$1/2$ is the unique maximum point if and only if $\lambda_c(]x,1-x[)>0$
for all $x<1/2$. This proves the last point.
\end{proof}

\subsection{Intricacy as a function of the profile}

Let us set
\[
\Gamma:=\left\{ h:[0,1]\mapsto [0,1]: \
h(0)=0,\; t\mapsto h(t) \text{ is non-decreasing and $1$-Lipschitz}\right\}
\]
and for any real number $x\in[0,1]$
\[
\Gamma_{x}:=\left\{ h\in \Gamma: h(1)=x\right\}.
\]
These sets are endowed with the partial order: $h\leq g$
if and only if $h(t)\leq g(t)$ for all $t\in[0,1]$. Each $\Gamma_x$
has a unique maximal element: the previously introduced
ideal entropy profile, $h^*_x(t)=t\wedge x$.

\begin{lemma}\label{gamma}
For any $X\in\cX(d,N)$, the entropy profile $h_X$, defined according
to Def. \ref{def:profile}, belongs to $\Gamma$.
\end{lemma}
\begin{proof}
Let $X\in X(d,N)$. Setting $I:=\{1,\dots,N\}$ and
\[
H_k:= \frac{1}{\binom{N}{k}}
     \sum_{S\subset I,\, |S|=k}       \H(X_S), \qquad k=0,1,\dots,N
\]
we must prove that
 $$
    0 = H_0\leq H_1 \leq\dots \leq H_{N}=\H(X), \qquad
    H_{k+1}-H_k\leq \log d, \quad 0\leq k<N.
 $$
The equalities $H_0=0$ and $H_N=\H(X)$ are obvious.
Let $0\leq k<N$ and compute:
\[
 \begin{split}
    H_{k+1} & = \frac{1}{\binom{N}{k+1}} \sum_{|S|=k+1}
            \H(X_{S})= \frac{1}{\binom{N}{k+1}} \sum_{|S|=k}
            \frac{1}{{k+1}} \sum_{i\in S^c} \H(X_{S\cup\{i\}})  \\
            & \leq \frac{k! (N-k-1)!}{N!}
                 \sum_{|S|=k}  (N-k) (\H(X_S)+\log d)  = H_k+\log d
 \end{split}
 \]
where $\H(X_{S\cup\{i\}})\leq \H(X_S)+\H(X_i)$ by \eqref{leq} and $\H(X_i)\leq \log d$ by \eqref{orly}.
The same computation, since $\H(X_{S\cup\{i\}})\geq \H(X_S)$ by \eqref{geq}, proves $H_k\leq H_{k+1}$.
\end{proof}

Let for any $h\in\Gamma$
 $$
   G_N^c(h):=2\sum_{k=0}^N c^N_k\binom{N}{k} \, h(k/N)-h(1) = 2\,
   \E\left(h\left(\beta_N\right)\right)-h(1).
 $$
\begin{lemma}\label{optipro}
Fix $x\in[0,1]$.
\begin{enumerate}
 \item For all $X\in\cX(d,N)$
$
G^c_N(h_X) = \frac{\I(X)}{N\log d}.
$
 \item $h^*_x$ is the unique maximizer of $G^c_N$ in
$\Gamma_x$ and
$
G^c_N(h_x^*)=i^c_N(x).
$
 \item For arbitrary $h\in\Gamma_x$, we have
\begin{equation}\label{mutte}
\|h-h^*_x\|_{c,N} = |G^c_N(h)-G^c_N(h^*_x)|.
\end{equation}
\end{enumerate}
\end{lemma}

\begin{proof}
Since $\MI(X,Y)=\H(X)+\H(Y)-\H(X,Y)$, $c^I_S=c^I_{S^c}$, and $\sum_S c^I_S=1$, we obtain
 $$
   \I^c(X) =
 2\sum_{k=0}^N c^N_k\sum_{|S|=k}\H(X_S) - \H(X).
 $$
Hence, the intricacy can be computed from the entropy profile:
\begin{equation}\label{permu}
\frac{\I^c(X)}{N\log d} =
 2\sum_{k=0}^N c^N_k\binom{N}{k} h_X(k/n) - h_X(1)
= G^c_N(h_X)
\end{equation}
and (1) is proved.

A direct computation yields for arbitrary $h\in\Gamma_x$:
\[
 \begin{split}
    |G^c_N(h^*_x)-G^c_N(h)| & = G^c_N(h^*_x)-G^c_N(h)
       = 2\sum_{k=0}^N c^N_k\binom{N}{k}(h_x^*(k/N)-h(k/N)) \\ &=
   2\sum_{k=0}^N c^N_k\binom{N}{k}|h_x^*(k/N)-h(k/N)| = \|h^*_x-h\|_{c,N},
   \end{split}
\]
since each term is non-negative. This proves (3).

Observe that $G_N^c:\Gamma_x\to\R$ is monotone non-decreasing. Hence,
setting $x=\frac{\H(X)}{N\log d}$ and recalling \eqref{eq:max-intric2}
 $$
   \I^c(X)=G^c_N(h_X)\leq \sup_{h\in\Gamma_{x}} G^c_N(h)
      = G^c_N(h_x^*)=2\E\left(x\wedge \beta_N\right)-x=i^c_N(x).
 $$
Moreover, $\I^c$ being non-null, all $c^I_S$ are positive, $G^c_N$ is increasing and
$h^*_x$ is a maximizer. Uniqueness of the maximizer in $\Gamma_x$ follows from
\eqref{mutte}, and point (2) is proved.
\end{proof}

\subsection{Proof of Proposition \ref{violi}}

Formula \eqref{mar} follows from Lemma \ref{optipro}, since for $\frac{\H(X)}{N\log d}=x$
\[
\frac{\I^c(X)}{N\log d}=G_N^c(h_X)=G_N^c(h_x^*)+G_N^c(h_X)-G_N^c(h_x^*)=
i^c_N(x)-\|h_X-h^*_x\|_{c,N}.
\]
To prove \eqref{mari}, it is enough to use \eqref{mar},
together with the continuity of $i^c$ and the uniform convergence of
$i^c_N\to i^c$. Proposition \ref{violi} is proved.

\subsection{No system with the ideal profile}
We turn to the problem of maximizing $\I^c$ over $\cX(d,N)$ at fixed $N$ for a prescribed
value of the entropy $\H(X)$.
The above results show that a system $X\in\cX(d,N)$ such that
$h_X(k/N)=h^*_x(k/N)$ ($k=0,1,\dots,N$) with $x=\frac{\H(X)}{N\log d}$
would be an exact maximizer.
However, the next Lemma
shows that such $X$ cannot exist except if $K$ or $N-K$ are bounded, independently of $N$.
Thus, all we can hope is to find systems
which approach the ideal profile. This will be done in section \ref{appri}.
\begin{lemma}\label{lem:no-ideal-pro}
For each $d\geq2$, there exists $H_*=H_*(d)<\infty$ with the following property.
If $N\geq1$ and $Y_1,\ldots,Y_N$ are random variables taking values in $\{0,\ldots,d-1\}$
and defined on the same probability space such that, for some real number $H\in[0,N]$,
\[
\frac{\H(Y_{\sigma(1)},\ldots,Y_{\sigma(k)})}{\log d} = k\wedge H, \qquad \forall
\ \sigma\in\cS_N, \ \forall k=1,\ldots,N,
\]
then $H$ or $N-H\leq H_*$.
\end{lemma}
\begin{proof}
Let $K:=\lfloor H\rfloor$ and $\tilde K:=\lceil H\rceil$.
Without loss of generality, we assume that $K\geq 3$ and
we proceed by contradiction.
Let us condition on the variables $(X_3,\dots,X_{\tilde K})$ (in the following paragraphs
we simply write "conditional" for "conditional on $(X_3,\dots,X_{\tilde K})$). By assumption:
 \begin{itemize}
  \item $(X_1,X_2)$ belongs to $Z:=\{0,\dots,d-1\}^2$;
  \item each $X_i$, $\tilde K<i\leq N$, is a function of $X_1,X_2$ as the conditional
entropy of $(X_1,X_2,X_i)$ is not bigger than that of $(X_1,X_2)$. Moreover,
the conditional entropy of $X_i$ is $\log d$. Hence, each such $X_i$
  defines a partition $\mathcal Z_i$ of $Z$ into $d$ subsets.
  \item For any pair $i\ne j$ in $\{1,2,\tilde K+1,\dots,N\}$, $(X_i,X_j)$ has conditional entropy $(H-\tilde K+2)\log d$, strictly greater than that of $X_i$ or $X_j$, both equal to $\log d$.
In particular, $\mathcal Z_i\ne\mathcal Z_j$.
 \end{itemize}
Thus, we have an injection from $\{1,2,\tilde K+1,\dots,N\}$ into the set of partitions of
$Z$ into $d$ subsets. This implies:
 $$
      N-\tilde K+2 \leq \binom{d^2+d-1}{d-1}.
 $$
Thus $N-H\leq H_*(d):= \binom{d^2+d-1}{d-1}$.
\end{proof}

\section{Random Construction of approximate maximizers}
\label{appri}

Motivated by \eqref{mari}, we introduce the following
\begin{definition}\label{maxi}
Let $\I^c$ be some intricacy and let $x\in[0,1]$ and $d\geq2$.
\begin{enumerate}
\item The {\bf entropy-intricacy function} $\I^c(d,x)$ is:
\begin{equation}\label{I(d,x)}
\sup\left\{ \limsup_{N\to+\infty}\frac{\I^c(\mu_N)}{N\log d}: \ {\mu_N\in\cM(d,N)},
\ \lim_{N\to+\infty}\frac{\H(\mu_N)}{N\log d}=x \right\}.
\end{equation}
\item A sequence of systems $X^N\in\cX(d,N)$, $N\geq2$, is an
{\bf approximate $x$-maximizer} for $\I^c$ if
\[
\lim_{N\to\infty}\frac{\H(X^N)}{N\log d}=x \qquad {\rm and} \qquad
     \lim_{N\to\infty} \frac{\I^c(X^N)}{N\log d} = \I^c(d,x).
\]
\item $(X^N)_N$ is an {\bf approximate maximizer} for $\I^c$
if
\[
     \lim_{N\to\infty} \frac{\I^c(X^N)}{N\log d} = \max_{x\in[0,1]} \I^c(d,x).
\]
\end{enumerate}
\end{definition}

Proposition \ref{violi} established that $\I^c(d,x)\leq i^c(x)$.
Proposition \ref{duque} in this section shows that this inequality is in
fact an equality.

In the rest of this section we construct approximate $x$-maximizers
by choosing uniform distributions on {\it random supports} with the appropriate size:
since $\frac{\H(\mu^N)}{N\log d}$ must be close to $x$ and $\mu^N$
is uniform, then the size of the (random) support of $\mu^N$ must be close to $d^x$, see \eqref{orly}.
It turns out that this simple construction yields the desired results.

\begin{rem}
{\rm
In \cite{buza1} we have given a different definition of $\I^c(d,x)$, namely
\[
  \I^c(d,x):= \lim_{N\to\infty}
            \frac1{N\log d} \, \sup \left\{ \I^c(X): X\in\cX(d,N), \
            \left|\frac{\H(X)}{N\log d} - x\right| \leq \delta_N \right\}
\]
for any sequence $(\delta_N)_{N}$ of non-negative numbers converging to zero.
It is easy to see that this definition and \eqref{I(d,x)} actually coincide.
}
\end{rem}

\subsection{Sparse random configurations}

Let $N\geq 2$ and $0\leq M\leq N$ be integers. We denote
\[
\Lambda_{d,n}:=\{0,\ldots,d-1\}^n, \qquad \forall \ n\geq 1.
\]
We consider a family
$(W_i)_{i\in\Lambda_{d,M}}$ of i.i.d. variables, each uniformly distributed
on $\Lambda_{d,N}$, defined on a probability space
$(\Omega,\cF,\bbP)$. We define a {\it random} probability measure
on $\Lambda_{d,N}$
\begin{equation}\label{muN}
\mu^{N,M}(x) := d^{-M} \sum_{i\in\Lambda_{d,M}} \un{(x=W_i)}, \qquad x\in \Lambda_{d,N}.
\end{equation}
In what follows we consider random variables $X^{N,M}$ on $(\Omega,\cF,\bbP)$ such that
\begin{equation}\label{XN}
\bbP\left( \left. X^{N,M} = x \, \right| \, (W_i)_{i\in\Lambda_{d,M}} \right) = \mu^{N,M}(x),
\qquad x\in \Lambda_{d,N}.
\end{equation}
In other words,
\[
{\it conditionally \ on \ } (W_i)_{i\in\Lambda_{d,M}}, \quad X^{N,M} \ {\it has \ law } \ \mu^{N,M}.
\]
\noindent
We are going to prove the following
\begin{prop}\label{duque}
For integers $N\geq1$, $0\leq M\leq N$, let $X^{M,N}$ be the random systems defined above. Let $x\in [0,1]$.
For any intricacy $\I^c$ we have, a.s. and in $L^1$
\begin{equation}\label{guo2}
\lim_{N\to+\infty} \frac{\I^c(\mu^{N,\lfloor xN\rfloor })}{N\log d}
= i^c(x)
\end{equation}
and
\begin{equation}\label{guo}
\lim_{N\to+\infty} \frac{\H(\mu^{N,\lfloor xN\rfloor })}{N\log d} = x.
\end{equation}
\end{prop}

\begin{rem}\label{stress1}{\rm
We stress that in the following we denote
\begin{equation}\label{abuse}
\I^c(X^{N,M})=\I^c(\mu^{N,M}), \qquad \H(X^{N,M})=\H(\mu^{N,M})
\end{equation}
and that all these expressions are random variables which depend
on $(W_i)_{i\in\Lambda_{d,M}}$. In other words, \eqref{abuse} indicates entropy
and intricacy of the law of $X^{N,M}$ {\it conditionally on }
$(W_i)_{i\in\Lambda_{d,M}}$. This abuse of notation seems necessary, to keep notation reasonably readable.
See also Remark \ref{boune} below.
}
\end{rem}

\subsection{Average intricacy of sparse random configurations}
We recall that $N\geq 2$, $M$ is an integer between $1$ and $N$
and  $\cS_N$ denotes the set of permutations of $\{1,\dots,N\}$.
By Lemma \ref{optipro},
$\I^c(X^{N,M})=2\sum_{k=0}^N c^N_k\binom{N}{k} h_X(k/N)-h_X(1)$, hence we get:
\begin{equation}\label{svedese}
\E\left({\I^c(X^{N,M})}\right)=
\frac1{N!}\sum_{\sigma\in\cS_N} 2\sum_{k=1}^N c^N_k \, \binom{N}{k} \,
\E\left(\H(X^{N,M}_{\{\sigma(1),\ldots,\sigma(k)\}})\right)
- \E\left(\H(X^{N,M})\right).
\end{equation}
We are going to simplify this expression by exploiting the symmetries of our
construction.
\begin{lemma}
The random vector $X^{N,M}=(X_1^{N,M},\ldots,X_N^{N,M})\in\cX(d,N)$ is exchangeable, i.e.
for all $\sigma\in\cS_N$ and any $\Phi:\Lambda_{d,N}\mapsto\R$
\[
\E\left(\Phi(X_{\sigma(1)}^{N,M},\ldots,X_{\sigma(N)}^{N,M})\right) =
\E\left(\Phi(X_1^{N,M},\ldots,X_N^{N,M})\right).
\]
\end{lemma}
\begin{proof}
Note that every $\sigma\in\cS_N$ induces a permutation $\Sigma_\sigma:\Lambda_{d,N}\mapsto\Lambda_{d,N}$
\[
\Sigma_\sigma(x_1,\ldots,x_N)=(x_{\sigma(1)},\ldots,x_{\sigma(N)}), \qquad x\in\Lambda_{d,N}.
\]
In particular, $(X_{\sigma(1)}^{N,M},\ldots,X_{\sigma(N)}^{N,M})=\Sigma_\sigma (X^{N,M})$ has,
conditionally on $(W_i)_{i\in\Lambda_{d,M}}$,
distribution
\[
\Sigma_\sigma^*\mu^{N,M}(x) := d^{-M} \sum_{i\in\Lambda_{d,M}} \un{(x=\Sigma_\sigma(W_i))}, \qquad x\in \Lambda_{d,N}.
\]
However, $(\Sigma_\sigma(W_i))_{i\in\Lambda_{d,M}}$ has the same distribution as $(W_i)_{i\in\Lambda_{d,M}}$.
Therefore we conclude.
\end{proof}
\begin{rem}\label{boune}
{\rm  Notice that $\Sigma_\sigma^*\mu^{N,M}$ has same {\it law} as $\mu^{N,M}$,
but in general the two measures are {\it not} a.s. equal. In other words,
$(X_1^{N,M},\ldots,X_n^{N,M})$ is exchangeable but {\it not} exchangeable conditionally on $(W_i)_{i\in\Lambda_{d,M}}$.
}
\end{rem}
\noindent
In particular, for all $k\in\{1,\ldots,N\}$ and $\sigma\in\cS_N$
\begin{equation}\label{roth}
\E\left(\H(X^{N,M}_{\{\sigma(1),\ldots,\sigma(k)\}})\right)=
\E\left(\H(X^{N,M}_{\{1,\ldots,k\}})\right),
\end{equation}
and we obtain by \eqref{svedese}
\begin{equation}\label{svedese2}
\E\left({\I^c(X^{N,M})}\right)=
2\sum_{k=1}^N c^N_k \, \binom{N}{k} \,
\E\left(\H(X^{N,M}_{\{1,\ldots,k\}})\right) - \E\left(\H(X^{N,M})\right).
\end{equation}
\begin{lemma}\label{acero}
Let $y\in \Lambda_{d,k}$, $k\in\{1,\ldots,N\}$ and set
\begin{equation}\label{nu}
\nu(y):=\sum_{z\in\Lambda_{d,N-k}} \mu^{N,M}(y,z).
\end{equation}
Then $d^M\cdot\nu(y)$ is a binomial variable with parameters $(d^M,d^{-k})$.
\end{lemma}
\begin{proof}
Notice that, conditionally on $(W_i)_{i\in\Lambda_{d,M}}$, $X^{N,M}_{\{1,\ldots,k\}}=(X^{N,M}_1,\ldots,X^{N,M}_k)\in\Lambda_{d,k}$
has distribution
\[
\bbP\left(\left.X^{N,M}_{\{1,\ldots,k\}}=y \, \right| \, (W_i)_{i\in\Lambda_{d,M}}\right) = \sum_{z\in\Lambda_{d,N-k}}
\mu^{N,M}(y,z) = d^{-M} \sum_{i\in\Lambda_{d,M}} \sum_{z\in\Lambda_{d,N-k}} \un{((y,z)=W_i)},
\]
where $y\in\Lambda_{d,k}$ and
$(y,z)\in \Lambda_{d,k}\times\Lambda_{d,N-k}=\Lambda_{d,N}$.
For fixed $y\in\Lambda_{d,k}$, the family
\[
T_i:=\sum_{z\in\Lambda_{d,N-k}} \un{((y,z)=W_i)}, \qquad i\in\Lambda_{d,M}
\]
is an i.i.d. family of Bernoulli variables with parameter $d^{-k}$. Indeed, if
$\Pi_{N\mapsto k}:\Lambda_{d,N}\mapsto\Lambda_{d,k}$ is the natural projection, then the law
of $\Pi_{N\mapsto k}(W_i)$ is uniform on $\Lambda_{d,k}$, so that
\[
\bbP(T_i=1)=\bbP(\Pi_{N\mapsto k}(W_i)=y)=d^{-k}.
\]
Hence, for all $y\in\Lambda_{d,k}$, $d^M\cdot\nu(y)=\sum_{i\in\Lambda_{d,M}}T_i$
is the sum of $d^M$ independent Bernoulli variables with parameter $d^{-k}$,
i.e. a binomial variable with parameters $(d^M,d^{-k})$.
\end{proof}
\noindent
Let us denote from now on by $B_k$ a binomial variable with parameters $(d^M,d^{-k})$.
Set
\[
\varphi(x):=-\frac{x\, \log x}{\log d}, \quad \forall x>0, \qquad \varphi(0):=0.
\]
Notice that the function $\psi(x):=-(1+x)\log(1+x) + x + \frac{x^2}2$ satisfies
\[
\psi(0)=\psi'(0)=0, \quad \psi''(x)\geq 0, \quad \forall x\geq 0,
\]
so that $\psi(x)\geq 0$ for $x\geq 0$. Moreover, $\varphi(1+x)\geq 0$ if $x\in [-1,0]$.
Hence, for all $x\geq-1$,
 \begin{equation}\label{eq:phi-bound}
  \varphi(1+x)\geq -\un{(x>0)}\frac{x+x^2/2}{\log d}.
 \end{equation}

Now, by \eqref{nu}
\[
\H\left(X^{N,M}_{\{1,\ldots,k\}}\right) = -\sum_{y\in\Lambda_{k,N}} \nu(y)\,
\log\nu(y),
\]
then we obtain by Lemma \ref{acero} that
\begin{equation}\label{orchidea}
h_k:=\frac{1}{\log d}\, \E\left(\H\left(X^{N,M}_{\{1,\ldots,k\}}\right)\right) =
d^k \, \E\left(\varphi\left(B_k \, d^{-M}\right)\right).
\end{equation}
\begin{lemma}\label{orly1}
We have, for any $0\leq k\leq M$,
\[
h_k = k+\E\left(\varphi(B_k\, d^{k-M})\right)= M +d^{k-M} \,
\E\left(\varphi(B_k)\right).
\]
\end{lemma}
\begin{proof}
These identities follow from the formulae:
$\E(B_k)=d^{M-k}$, $\varphi(d^{-j})=jd^{-j}$ and
$\varphi(\alpha x)=\alpha\varphi(x)+x\varphi(\alpha)$ for $\alpha>0$
applied to $\varphi(B_kd^{k-M}\cdot d^{-k})$ and
$\varphi(B_k\cdot d^{-M})$.
\end{proof}
\begin{lemma}\label{flora}
Away from $k=M$, the entropy is nearly constant:
\[
k-2d^{\frac{k-M}2}\leq h_k\leq k,
\qquad k=1,\ldots,M,
\]
\[
M-d^{M-k}\leq h_k\leq M,
\qquad k=M+1,\ldots,N.
\]
\end{lemma}
\begin{proof}
The upper bounds are easy. Indeed, for $k\leq M$ one uses \eqref{leq}, while
for $k>M$ we notice that the support of $\mu^{N,M}$ has cardinality at most
$d^M$, and apply \eqref{orly} to conclude.

Recall that $B_k$ is binomial with parameters $(d^M,d^{-k})$. Then
$\E(B_k)=d^{M-k}$ and ${\rm Var}(B_k)=d^Md^{-k}(1-d^{-k})$. If
we define $J_k:=B_k\cdot d^{k-M}-1$ then we obtain
\[
\E\left(J_k^2\right) = d^{2(k-M)} \, {\rm Var}(B_k)=
d^{k-M}-d^{-M}\leq d^{k-M}.
\]
Hence, using \eqref{eq:phi-bound} we get
\[
\begin{split}
\E\left(\varphi(B_k\cdot d^{k-M})\right) &  = \E\left(\varphi(1+J_k)\right)
\geq -\frac1{\log d}\,
\E\left(\un{(J_k>0)}\left(J_k+\frac{J_k^2}2\right)\right)
\\ & \geq -\E\left(|J_k|+J_k^2\right)
 \geq -\left(\sqrt{\E\left(J_k^2\right)}+\E\left(J_k^2\right)\right)\geq-
2d^{\frac{k-M}2},
\end{split}
\]
since $\E(|J_k|)\leq \sqrt{\E(J_k^2)}$ by Cauchy-Schwartz and $d^{\frac{k-M}2}\leq d^{k-M}\leq 1$.
By Lemma \ref{orly1} we obtain the desired lower bound for $k\leq M$.

Let us consider now the regime $k> M$.
We have
\[
h_k=d^k \, \E(\varphi(B_k\, d^{-M})) = M +d^{k-M} \,
\E\left(\varphi(B_k)\right).
\]
If $B_k\in\{0,1\}$ then $\varphi(B_k)=0$. Note that \eqref{eq:phi-bound} implies that
$\varphi(B_k)\geq - (B_k-1)^2-(B_k-1)$ as the right hand side is zero whenever $B_k=0,1$
and less than $-(B_k-1)^2/2-(B_k-1)$ otherwise. Thus,
\[
\begin{split}
d^{k-M}\E\left(\varphi(B_k)\right) & \geq -d^{k-M}\E\left((B_k-1)^2+(B_k-1)\right)
=-d^{k-M}\E\left(B_k^2-B_k\right)
\\ & = -d^{k-M}\left(d^{M-k}+d^{2(M-k)}-d^{M-2k}-d^{M-k}\right)
\\ & \geq -d^{M-k}+ d^{-k} \geq -d^{M-k}.
\end{split}
\]
By Lemma \ref{orly1} we obtain the lower bound for $k>M$.
\end{proof}

\subsection{Estimation of the expected Intricacy}

\begin{lemma}\label{fili}
Let $x\in\, ]0,1[$, $M:=\lfloor xN\rfloor$ and $\alpha:=d^{1/2}>1$.
For all $N\geq 2$, using
the notation \eqref{lava0},
\begin{equation}\label{svedese3}
-2d^{-(N-M)}\leq \E\left(2\, (D_N\wedge M)-M\right)-
\frac{\E\left(\I^c(X^{N,M})\right)}{\log d} \leq
\E\left(4\alpha^{-|D_N-M|}\right).
\end{equation}
\end{lemma}
\begin{proof}
By \eqref{svedese2}, \eqref{orchidea} and \eqref{lava},
\[
\frac{\E\left({\I^c(X^{N,M})}\right)}{\log d}=
2\sum_{k=1}^N c^N_k \, \binom{N}{k} \,
h_k - h_N=2\, \E(h_{D_N})-h_N.
\]
So,
$$
 \frac{\E\left({\I^c(X^{N,M})}\right)}{\log d}-\E(2(D_N\wedge M)-M)
 = 2\E(h_{D_N}-D_N\wedge M)+M-h_N
$$
We conclude by Lemma \ref{flora}.
\end{proof}

\begin{lemma}\label{bled}
Let $x\in\, ]0,1[$ and $M:=\lfloor xN\rfloor$. Then for
$\alpha:=d^{1/2}>1$ and some constant $C>0$ we have
\begin{equation}\label{svedese4}
\E\left(\alpha^{-|D_N-M|}\right)\leq \frac C{\sqrt N}, \qquad
\forall \, N\geq 1.
\end{equation}
\end{lemma}
\begin{proof}
To ease notation, in this proof we drop the subscript
$c$ from $W_c$.
By \eqref{mainc} and \eqref{lava}, we have that
\begin{equation}\label{pica2}
\bbP(D_N=k) = \binom Nk \, c^N_k
=\binom Nk \, \E\left(W^{k} \, (1-W)^{N-k} \right), \qquad k=0,\ldots,N.
\end{equation}
We claim that for all $1\leq k\leq \lfloor \frac N2\rfloor$ we have
$$
\bbP(D_N=k-1) \leq \bbP(D_N=k).
$$
Indeed,
\[
\begin{split}
& \bbP(D_N=k) - \bbP(D_N=k-1)  =
\\ & = \frac{N!}{k! \, (N-k+1)!} \,
\E\left(W^{k-1} \, (1-W)^{N-k} \left[ (N+1-k)W-k(1-W)\right]\right).
\\ & = \frac{N!}{k! \, (N-k+1)!} \,
\E\left(W^{k-1} \, (1-W)^{N-k} \left[ (N+1)W-k\right]\right).
\end{split}
\]
By the symmetry of $\lambda_c$ w.r.t. $1/2$ we have, since $2k\leq N+1$,
\[
\begin{split}
& \E\left(W^{k-1} \, (1-W)^{N-k} \left[ (N+1)W-k\right]\right)
= 2^{-N+1} \, \left(\frac{N+1}2-k\right) \bbP(W=1/2) + \\ & +
\E\left(W^{k-1} \, (1-W)^{N-k} \left[ (N+1)W-k\right]
\un{(W>1/2)}\right)+
\\ & +\E\left(W^{N-k} \, (1-W)^{k-1} \left[ (N+1)(1-W)-k\right]
\un{(W>1/2)}\right)
\\ & \geq  \E\Big( \left(W(1-W)\right)^{k-1}
\un{(W>1/2)} \cdot
\\ & \qquad \cdot  \left[
 \left( (N+1)W-k\right)W^{N-2k+1}+
 \left( (N+1)(1-W)-k\right)(1-W)^{N-2k+1}
\right]\Big).
\end{split}
\]
This expectation is nonnegative. Indeed, the first term is nonnegative
as $2k\leq N+1$ and $W>1/2$. We assume the second term to be negative,
otherwise we are done. As $1-W\leq W$, we get:
\[
\begin{split}
& \bbP(D_N=k) - \bbP(D_N=k-1)  = \\ & =
\E\left[ \un{W>1/2} ((N+1)W-k+(N+1)(1-W)-k)(1-W)^{N-2k+1}\right]
\\ & \geq \E\left[ \un{W>1/2} ((N+1)-2k)W^{N-2k+1}\right]\geq 0,
\end{split}
\]
proving the claim.

We set $L:=\lfloor \frac N2\rfloor$.
Then $N-L\geq L\geq \frac N2-1$ and
by \eqref{pica2} we obtain for all $0\leq k\leq N$
\[
\begin{split}
\binom Nk \, c^N_k & \leq \binom N{L} \, c^N_{L} =
\frac{N!}{L!(N-L)!} \, \E\left(W^{L}\, (1-W)^{N-L}\right)
\\ & \leq \frac{N!}{L!(N-L)!} \,
\E\left(W^{L}\, (1-W)^{L}\right) \leq  \frac{N!}{L!(N-L)!} \,
4^{-L} \leq \frac{N!}{L!(N-L)!} \, 2^{-N+2}.
\end{split}
\]
By Stirling's formula $n!=\sqrt{2\pi n}(n/e)^n(1+\mathcal O(1/n))$,
there is a constant $C\geq 0$ such that
\[
2^{-N} \, \frac{N!}{L!(N-L)!} \leq C \, N^{-\frac12}, \qquad N\to+\infty, \
L=\left\lfloor \frac N2\right\rfloor.
\]
Then, we obtain for some constants $C_1,C_2\geq 0$
\[
\begin{split}
\E\left(\alpha^{-|D_N-M|}\right)& = \sum_{k=0}^N \binom Nk \, c^N_k \,
\alpha^{-|k-M|} \leq C_1 \, N^{-\frac12} \sum_{k=0}^N \alpha^{-|k-M|}
\\ & \leq C_1 \, N^{-\frac12} \, 2\sum_{k=0}^{+\infty} \alpha^{-k}
= C_2 \, N^{-\frac12},
\end{split}
\]
and the proof is finished.
\end{proof}

\begin{proof}[Proof of Proposition \ref{duque}]
Let $x\in\, ]0,1[$ and $M:=\lfloor xN\rfloor\geq 1$ ($N$ is large).
By Lemma \ref{flora} for $k=N$
\[
\begin{split}
\E\left(\left|\frac MN-\frac{\H(X^{N,M})}{N\log d}\right|\right) & =
\E\left(\frac MN-\frac{\H(X^{N,M})}{N\log d}\right)= \E\left(\frac{M-h_N}N\right)
\\ & \leq \frac{d^{M-N}}N
\leq d^{-N(1-x)}.
\end{split}
\]
Thus,
\begin{equation}\label{dern}
\sum_{N\geq 1} \E\left(\left|\frac MN-\frac{\H(X^{N,M})}{N\log d}\right|\right) <+\infty,
\end{equation}
therefore a.s.
\[
\sum_{N\geq 1} \left|\frac MN-\frac{\H(X^{N,M})}{N\log d}\right| <+\infty,
\]
and in particular a.s.
\[
\lim_{N\to+\infty} \left|\frac MN-\frac{\H(X^{N,M})}{N\log d}\right| =0.
\]
Setting $x_N:=\frac{\H(X^{N,M})}{N\log d}$, we have obtained
\[
\lim_{N\to+\infty} x_N= \lim_{N\to+\infty} \frac{\lfloor xN\rfloor}N =x
\]
a.s. and in $L^1$, namely we have proven \eqref{guo}.
Now, let us set observe that,
by Proposition \ref{violi}, this gives
 $$
   |i^c_N(x_N)-i^c(x)|\leq N^{-1/2}+|x_N-x|\to 0
 $$
again a.s. and in $L^1$.
On the other hand, by \eqref{eq:max-intric2} and
by Lemmas \ref{fili}, \ref{bled}
 $$
  \sum_{N\geq1} \E\left(\left| i^c_N(x_N)-\frac{\I^c(X^{N,M})}{N\log d} \right|\right)
  \leq C\sum_{N\geq1} N^{-3/2} < \infty.
 $$
Arguing as above,
it follows that $\frac{\I^c(X^{N,M})}{N\log d}\to i^c(x)$ a.s. and
in $L^1$.
This proves \eqref{guo2} and
concludes the proof of Proposition \ref{duque}.
\end{proof}

\section{Results for Arbitrary Intricacies}

We now collect our results to state the generalizations of
Theorems \ref{thm:max} and \ref{thm:profile} for arbitrary intricacies.
We consider some non-null intricacy $\I^c$. Let $\lambda_c$ be the
associated probability measure on $[0,1]$ according to
Proposition \ref{pro:special}.
Recall from Def. \ref{maxi} that the corresponding
entropy-intricacy function $\I^c(d,x)$ is:
 $$
     \I^c(d,x):=\sup\left\{\limsup_{N\to\infty} \frac{\I^c(X^N)}{N\log d}: X^N\in\cX(d,N)\text{ s.t. }
          \lim_{N\to\infty} \frac{\H(X^N)}{N\log d}=x\right\}.
 $$
We also recall that $i^c(x)$ and $i^c_N(x)$ have been defined in
eq. \eqref{eq:max-intric2} and \eqref{eq:max-intric}.

\begin{theorem}\label{thm:main}
$ $
\begin{enumerate}
 \item For any $N\geq1$, $X\in\cX(d,N)$,
 \[
 \frac{\I^c(X)}{N\log d} \leq i^c_N(x)
:=2\sum_{k=0}^N c^N_k \, \binom{N}{k} \,
  (k/N)\wedge x -x.
  \]
 \item $\I^c(d,x)=i^c(x)=2\int_0^1 x\wedge t\ \lambda_c(dt)-x$.
The function $i^c$ is Lipschitz with constant $1$, concave
and symmetric: $i^c(1-x)=i^c(x)$.
 \item $i^c(1/2)=\max_{x\in[0,1]} i^c(x)$. Moreover $1/2$ is the unique maximum if and only
 if $1/2$ belongs to the support of $\lambda_c$.
 \item $|i^c(x)-i^c_N(x)|\leq N^{-1/2}$ for all $x\in[0,1]$.
\end{enumerate}
\end{theorem}
\noindent
Theorem \ref{thm:main} immediately follows from Propositions \ref{violi} and \ref{duque}.
We now consider the convergence of entropy profiles to the ideal profiles
for approximate maximizers, i.e.,
the generalization of Theorem \ref{thm:profile}:

\begin{theorem}\label{barry}
Let $\I^c$ be an intricacy.
\begin{enumerate}
 \item If $1/2$ is in the support of $\lambda_c$,
then any approximate maximizer $(X^N)_N$ for $\I^c$ satisfies:
\begin{equation}\label{marian}
\lim_{N\to+\infty}\frac{\H(X^N)}{N\log d} = \frac12, \qquad
\lim_{N\to+\infty}\sup_{\operatorname{supp}(\lambda_c)}
    |h_{X^N}-h^*_{1/2}| = 0.
\end{equation}
 \item Let $x\in[0,1]$ and let $(X^N)_N$ be an
approximate $x$-maximizer for $\I^c$. Then:
 \begin{equation}\label{eq:right-profile}
    \lim_{N\to\infty} \sup_{\operatorname{supp}(\lambda_c)}
    |h_{X^N}-h^*_x| = 0.
 \end{equation}
In particular, if $x\in\operatorname{supp}(\lambda_c)$ then
 \begin{equation}\label{eq:right-profile2}
    \lim_{N\to\infty} \sup_{[0,1]}
    |h_{X^N}-h^*_x| = 0.
 \end{equation}
\item If $x\in\operatorname{supp}(\lambda_c)$ then an
approximate $x$-maximizer $(X^N)_N$ for $\I^c$ is an
approximate $x$-maximizer for any other intricacy $\I^{c'}$.
\end{enumerate}
\end{theorem}
\begin{rem}{\rm
The extra assumption about the support of $\lambda_c$  cannot be dropped. Indeed,
for point (1) of Theorem, observe that, in the case of the $p$-symmetric intricacy,  $1/2$ is not in the
support of $\frac12(\delta_p+\delta_{1-p})$ for $p\ne1/2$ and
approximate maximizers of $\I^p$ satisfy only
\[
p\leq \liminf_N \frac{\H(X^N)}{N\log d} \leq  \limsup_N \frac{\H(X^N)}{N\log d} \leq 1-p.
\]
The entropy may accumulate on any point on the
interval $[p,1-p]$.

Notice however that for many intricacies, including the neural complexiy,
the support of $\lambda_c$ is the whole interval, making this assumption
satisfied for all $x\in[0,1]$.
}
\end{rem}

\begin{rem}
{\rm
In the setting of point (2) of Theorem \ref{barry},
if $x$ does not belong to the support of $\lambda_c$, then
one can prove with similar arguments that
 \begin{equation}\label{eq:right-profile3}
    \lim_{N\to\infty} \sup_{[0,a]\cup [b,1]}
    |h_{X^N}-h^*_x| = 0
 \end{equation}
where $a:=\sup([0,x]\cap{\rm supp}(\lambda_c))$,
$b:=\inf([x,1]\cap{\rm supp}(\lambda_c))$, with the
convention $\sup\emptyset:=0$ and $\inf\emptyset:=1$.
}
\end{rem}

\begin{rem}\label{cano}{\rm
Let $X^{N,M}$ the random system constructed in \eqref{muN} and \eqref{XN},
with $M:=\lfloor xN\rfloor$, $x\in\,]0,1[$. For the Edelman-Sporns-Tononi intricacy $\I$
we have that the support of the associated probability measure $\lambda$ is $[0,1]$, since
it is the Lebesgue measure by Lemma 3.8 of \cite{buza1}. Since $(X^{N,M})_N$ is a.s.
an approximate $x$-maximizer for $\I$ by Proposition \ref{duque}, then by point
(3) of Theorem \ref{barry} a.s. this sequence is an approximate $x$-maximizer
{simultaneously} for all intricacies $\I^c$.

This has the following consequence for approximate maximizers (i.e., without
entropy constraints).
An approximate maximizer for some intricacy  $\I^c$ where $1/2\notin\operatorname{supp}(\lambda_c)$ is not necessarily an approximate maximizer for another intricacy. But
an approximate $1/2$-maximizer for any intricacy is automatically an approximate
maximizer for all intricacies.
}
\end{rem}

\begin{proof}[Proof of Theorem \ref{barry}]
Let us set for simplicity of notation:
 \begin{equation}\label{notation}
x_N:=\frac{\H(X^N)}{N\log d}, \qquad
\I_N:=\frac{\I^c(X^N)}{N\log d}.
 \end{equation}
If $1/2$ is in the support of $\lambda_c$, then it is
the unique point where $i^c(x)$ achieves its maximum.
Then, Theorem \ref{thm:main} implies that no $x\ne 1/2$ can
be an accumulation point of $x_N,N\geq1$. Thus an
approximate maximizer is an approximate $1/2$-maximizer.

It is now enough to prove point (2).
By definition of approximate $x$-maximizers, $x_N\to x$
and $\I_N\to i^c(x)$. Using point (2) of Proposition \ref{violi},
it follows that $|\I_N-i^c_N(x)|\to 0$ and by \eqref{mar} we have
$\|h_{X^N}-h_x^*\|_{c,N} \to 0$. Notice now that for any  $K$-Lipschitz
function $f:[0,1]\to\mathbb R$, by \eqref{oo}
 $$
    \left|\E(f(\beta_N)) - \E(f(W_c))\right| \leq \frac K{\sqrt N}
 $$
As entropy profiles are $1$-Lipschitz, we obtain
 $$
   \int |h_{X^N}-h_x^*| \, d\lambda_c \to 0.
 $$
As all functions $(h_{X^N}-h_x^*)_N$ are $2$-Lipschitz,
\eqref{eq:right-profile} follows by a routine argument.

Assuming now $x\in\operatorname{supp}(\lambda_c)$,
$\lim_{N\to\infty} h_{X^N}(x)=h_x^*(x)=x$. On
the one hand, as $h_{X^N}(0)=0$ and $h_{X^N}$ is $1$-Lipschitz,
it follows that the convergence $\lim_{N\to\infty} h_{X^N}(t)=h_x^*(t)=t$
occurs for all $t\in[0,x]$. On the other hand, all $h_{X^N}$ being non-decreasing,
$h_{X^N}(x)=x\leq h_{X^N}(t)
\leq h_{X^N}(1)\to x$. Hence the previous convergence occurs for all $x\in[0,1]$,
proving \eqref{eq:right-profile2}.

Finally, let us prove point (3). By \eqref{eq:right-profile2} the profiles $h_{X^N}$
converge to $h^*_x$ uniformly on $[0,1]$. Let $\I^{c'}$ be any other intricacy. By
uniform convergence we have $\|h_{X^N}-h_x^*\|_{c',N}\leq\sup_{[0,1]}|h_{X^N}-h_x^*|\to 0$.
By \eqref{mar} and Lemma \ref{opti-I}
\[
\frac{\I^{c'}(X)}{N\log d} = i^{c'}_N(x)-\|h_X-h_x^*\|_{c',N}\to i^{c'}(x), \qquad N\to+\infty.
\]
By Theorem \ref{thm:main}, $i^{c'}(x)=\I^{c'}(d,x)$ and therefore
$(X^N)_N$ is an approximate $x$-maximizer for $\I^{c'}$.
\end{proof}

We have
the following consequence for approximate $x$-maximizers.
We recall that $\H(Y\,|\,Z)$ denotes the
conditional entropy, see the Appendix.
\begin{theorem}\label{ott}
Suppose that $x\in{\rm supp}(\lambda_c)$ and let $(X^N)_N$ be an
approximate $x$-maximizer for some $0\leq x\leq 1$. Then
\begin{enumerate}
\item If $y\in\,]0, x]$ then for all $\varepsilon>0$
\begin{multline*}
\lim_{N\to+\infty} \frac{1}{\binom{N}{\lfloor yN\rfloor}}
  \#\biggl\{S\subset\{1,\ldots,N\}: |S|= \lfloor yN\rfloor, \\
     (1-\varepsilon)|S|\log d < \H(X_S^N) \leq |S|\log d\biggr\} =1.
\end{multline*}
\item If $y\in [x,1[$ then for all $\varepsilon>0$
\[
\lim_{N\to+\infty} \frac{1}{\binom{N}{\lfloor yN\rfloor}}
\#\{S\subset\{1,\ldots,N\}: |S|= \lfloor yN\rfloor, \
\H(X^N\, |\, X_S^N)<\varepsilon xN\log d\} =1.
\]
\end{enumerate}
\end{theorem}
\noindent
This result can be loosely interpreted as follows: as $N\to+\infty$,
\begin{enumerate}
\item if $y\in\,]0, x]$ then for {\it almost} all subsets $S$ with $|S|= \lfloor yN\rfloor$,
$X_S$ is {\it almost} uniform;
\item if $y\in [x,1[$ then for {\it almost} all subsets $S$ with $|S|= \lfloor yN\rfloor$,
$X$ is {\it almost} a function of $X_S$.
\end{enumerate}
This follows from the relation between entropy and
conditional entropy on one side and independence versus
dependence on the other side, see the Appendix.
\begin{proof}[Proof of Theorem \ref{ott}] Let $y\in\, ]0,1[$.
By \eqref{eq:right-profile2}, $h_{X^N}(y)\to h^*x(y)=x\wedge y$ as
$N\to+\infty$. By the definition \ref{def:profile} of $h_{X^N}$, we obtain,
setting $k_N:=\lfloor yN\rfloor$,
\[
\frac{1}{\binom{N}{k_N}}
     \sum_{|S|=k_N} \left|h^*_x(y)-\frac{\H(X_S)}{N\log d}\right| =
h^*_x(y) - \frac{1}{\binom{N}{k_N}}
     \sum_{|S|=k_N} \frac{\H(X_S)}{N\log d} = h^*_x(y) - h_{X^N}(y)\to 0,
\]
since all terms in the sum are non-negative by Lemma \ref{gamma}.
Let $\cZ_N$, defined on $(\Omega,\cF,\bbP)$, be a random subset
of $\{1,\ldots,N\}$ defined by
\[
\bbP(\cZ_N=S) = \frac{1}{\binom{N}{k_N}}, \qquad {\rm if} \
|S|=k_N.
\]
Then the above formula can be rewritten as follows
\[
\lim_{N\to+\infty}
\E\left(\left|h^*_x(y)-\frac{\H(X_{\cZ_N})}{N\log d}\right|\right) = 0.
\]
Since $L^1$ convergence implies convergence in probability, we
obtain
\[
\lim_{N\to+\infty}
\bbP\left(\left|h^*_x(y)-\frac{\H(X_{\cZ_N})}{N\log d}\right|>\varepsilon
h^*_x(y)\right)=0, \qquad \forall\, \varepsilon>0.
\]
This readily implies the Theorem, by recalling that
$\H(X^N\, |\, X_S^N)=\H(X^N)-\H(X_S^N)$ and that
$\frac{\H(X^N)}{N\log d}\to x$ by assumption.
\end{proof}

\appendix

\section{Entropy}

In this Appendix, we recall needed facts from basic information
theory. The main object is the entropy functional which may be
said to quantify the randomness of a random variable. We refer
to \cite{InfoTheory} for more background.

Let $X$ be a random variable taking values in a finite space $E$.
We define the {\it entropy} of $X$
\[
\H(X) := -\sum_{x\in E} P_X(x) \, \log(P_X(x)), \qquad
P_X(x):=\bbP(X=x),
\]
where we adopt the convention
\[
0\cdot \log(0) = 0\cdot \log(+\infty)=0.
\]
We recall that
\begin{equation}\label{orly}
0\leq \H(X) \leq \log |E|,
\end{equation}
More precisely, $\H(X)$ is minimal iff $X$ is a constant,
it is maximal iff $X$ is uniform over $E$. To prove \eqref{orly}, just notice that
since $\varphi\geq 0$ and $\varphi(x)=0$ if and only if $x\in\{0,1\}$, and
by strict convexity of $x\mapsto \varphi(x)=x\log x$ and Jensen's inequality
\[
\begin{split}
\log |E|-H(X) & = \frac1{|E|}\sum_{x\in E} P_X(x)\, |E| \left(\log(P_X(x))+ \log |E| \right)
\\ & = \frac1{|E|}\sum_{x\in E} \varphi\left(P_X(x)\, |E|\right) \geq \varphi
\left(\frac1{|E|}\sum_{x\in E} P_X(x)\, |E|\right) = \varphi(1)=0,
\end{split}
\]
with $\log |E|-H(X)=0$ if and only if $P_X(x)\, |E|$ is constant in $x\in E$.

\smallskip
If we have a $E$-valued random variable $X$ and a $F$-valued
random variable $Y$ defined on the same probability space,
with $E$ and $F$ finite, we can consider the vector
$(X,Y)$ as a $E\times F$-valued random variable The entropy of $(X,Y)$ is
then
\[
\H(X,Y) := -\sum_{x,y} P_{(X,Y)}(x,y) \, \log(P_{(X,Y)}(x,y)), \quad
P_{(X,Y)}(x,y):=\bbP(X=x,Y=y).
\]
This entropy $\H(X,Y)$ does not only depends on the (separate) laws
of $X$ and $Y$ but on the extent to which the "randomness of the
two variables is shared". The following notions formalize this
idea.

\subsection{Conditional Entropy}
The {\it conditional entropy} of $X$ given $Y$ is:
\[
\H(X\, |\, Y) := \H(X,Y) - \H(Y).
\]
We first claim that it is nonnegative.

Remark that $P_X(x)$ and $P_Y(y)$, defined in the obvious way,
are the marginal laws of $P_{(X,Y)}(x,y)$, i.e.
\[
P_X(x) = \sum_y P_{(X,Y)}(x,y), \qquad
P_Y(y) = \sum_x P_{(X,Y)}(x,y).
\]
In particular, $P_X(x)\geq P_{(X,Y)}(x,y)$ for all $x,y$. Therefore
\[
\sum_{x,y} P_{(X,Y)}(x,y) \, \log\left(\frac{P_{(X,Y)}(x,y)}{P_X(x)}\right)
\leq 0
\]
which yields
\[
\H(X,Y)=- \sum_{x,y} P_{(X,Y)}(x,y) \, \log P_{(X,Y)}(x,y)
\geq -\sum_x  P_X(x) \, \log P_X(x) = \H(X),
\]
i.e., $\H(X|Y)\geq0$, proving the claim. Therefore
\begin{equation}\label{geq}
\H(X,Y)\geq \max \{\H(X),\H(Y)\}.
\end{equation}
Moreover $\H(X,Y)=\H(X)$, i.e. $\H(Y|X)=0$, if and only if
$P_{(X,Y)}(x,y)=P_X(x)$ whenever
$P_{(X,Y)}(x,y)\ne 0$, which means that $Y$ is a function of $X$.

On the other hand,
 \begin{equation}\label{leq}
       \H(X,Y) \leq \H(X)+\H(Y)
 \end{equation}
with equality, i.e., $\H(Y|X)=\H(Y)$, if and only if $X$ and $Y$ are
independent. This can be shown by considering
the Kullback-Leibler divergence or relative entropy:
\[
I := \sum_{x,y} P_{(X,Y)}(x,y) \, \log\left(\frac{P_{(X,Y)}(x,y)}{P_X(x)\,P_Y(y)}\right).
\]
Since $\log(\cdot)$ is concave, by Jensen's inequality
\[
-I \leq \log \left( \sum_{x,y} P_{(X,Y)}(x,y) \,
\frac{P_X(x)\,P_Y(y)}{P_{(X,Y)}(x,y)}\right) =
\log \left( \sum_{x,y}P_X(x)\,P_Y(y)\right) = 0.
\]
By strict concavity, $I=0$ if and only if $P_{(X,Y)}(x,y)=P_X(x)\,P_Y(y)$
for all $x,y$, i.e., whenever $X$ and $Y$ are independent.

By the above considerations, $\H(X\, |\, Y)\in[0,\H(X)]$ is
a measure of the uncertainty associated with $X$ if
$Y$ is known. It is minimal iff $X$ is a function of $Y$ and
it maximal iff $X$ and $Y$ are independent.

\subsection{Mutual Information}
Finally, we recall the notion of
{\it mutual information} between two random variables $X$ and $Y$
defined on the same probability space:
\[
\begin{split}
\MI(X,Y) & := \H(X)+\H(Y)-\H(X,Y)
\\ &  = \H(X)-\H(X\, |\, Y) = \H(Y)-\H(Y\, |\, X)
\\ & = \sum_{x,y} P_{(X,Y)}(x,y) \, \log\left(\frac{P_{(X,Y)}(x,y)}{P_X(x)\,P_Y(y)}\right).
\end{split}
\]
This quantity is a measure of the common randomness of $X$ and $Y$.
By \eqref{geq} and \eqref{leq} we have
$\MI(X,Y)\in[0,\min\{\H(X),\H(Y)\}]$. $\MI(X,Y)$ is
minimal (zero) iff $X,Y$ are independent and maximal, i.e. equal
to $\min\{\H(X),\H(Y)\}$, iff one variable is a function of
the other.

\end{document}